\documentclass[12pt,reqno]{amsart}
\usepackage{amssymb}

 \rm

\renewcommand{\(}{\left(}
\renewcommand{\)}{\right)}
\renewcommand{\[}{\left[}
\renewcommand{\]}{\right]}

\newcommand{\x}{\times}
\renewcommand{\bar}{\overline}
\newcommand{\abs}[1]{\left\lvert#1\right\rvert}
\newcommand{\norm}[1]{\left\lVert#1\right\rVert}
\newcommand{\st}{\:|\:}

\renewcommand{\phi}{\varphi}

\renewcommand{\Re}{{\mathrm{Re}\,}}

\newcommand{\CC}{\mathbb{C}}
\newcommand{\RR}{\mathbb{R}}

\newcommand{\NN}{\mathbb{N}}

\newcommand{\M}{\mathcal{M}}
\newcommand{\F}{\mathcal{F}}
\newcommand{\G}{\mathcal{G}}

\newcommand{\I}{\mathrm{I}}
\newcommand{\II}{\mathrm{II}}
\newcommand{\Res}{\operatorname{Res}}

\theoremstyle{plain}
\newtheorem{thm}{Theorem}[section]

\newtheorem{ex}[thm]{Example}

\theoremstyle{definition}
\newtheorem{defn}[thm]{Definition}

\theoremstyle{remark}
\newtheorem{rem}[thm]{Remark}

\title{The Brylinski beta function of a double layer}
\author{Pooja~Rani and M.~K.~Vemuri}
\thanks{Pooja Rani is supported by DST-INSPIRE Fellowship/2019/IF190300}
\address{Department of Mathematical Sciences, IIT (BHU), Varanasi 221 005.}
\date{\today}
\date{1 March, 2023}

\begin{document}

\begin{abstract}
An analogue of Brylinski's knot beta function is defined for a
compactly supported (Schwartz) distribution $T$ on $d$-dimensional
Euclidean space.  This is a holomorphic function on a right
half-plane.  If $T$ is a (uniform) double-layer on a compact smooth
hypersurface, then the beta function has an analytic continuation to
the complex plane as a meromorphic function, and the residues are
integrals of invariants of the second fundamental form.  The first few
residues are computed when $d=2$ and $d=3$.
\end{abstract}

\maketitle


\section{Introduction}\label{S:intro}

In \cite{bry}, Brylinski introduced the beta function of a geometric knot in
$\RR^3$.  He was partly motivated by the desire to give a definition of
M\"obius energy (see \cite{fre,OHara1991241}) independent of an arbitrary
``renormalization''.  However, he also gave some beautiful formulae for the
first few residues of his beta function.  They turn out to be integrals of
polynomials in the curvature, torsion and their derivatives.

The basic idea behind Brylinski's beta function is much more general.  For
instance, in \cite{bbs}, Fuller and Vemuri defined and studied the beta
function of a compact submanifold of Euclidean space
(see also \cite{o2018regularized, err-o2018regularized}).  More generally, one
can define (for $\Re s >> 0$) the beta function of a compactly supported
(Schwartz) distribution on Euclidean space.  This function may not, in
general, have an interesting analytic continuation to the whole $s$-plane.
If the distribution is a (uniform) {\em single-layer}, i.e.
the surface measure of a compact submanifold, then the beta function
coincides with the one studied in \cite{bbs}, which does have an interesting
analytic continuation.  In this work, we study the beta-function of a
(uniform)
{\em double-layer}, i.e.\ a uniform distribution of unit-normal dipoles along
a hypersurface.

\section{The Brylinski beta function}

Let $T$ be a compactly supported distribution on $\RR^d$.  Since $T$
has finite order, there exists $k \in \NN$
such that $T$ extends to a continuous linear functional on $C^k(\RR^d)$.
Observe that if $\Re s > 2k$ and $N_s(u)=\norm{u}^s$ then
$N_s \in C^{2k}(\RR^d)$, and
so $N_s*T \in C^k(\RR^d)$.

\begin{defn}
The Brylinski beta function of $T$ is the function $B_T(s)$, defined for
$\Re s > 2k$ by
\begin{equation*}
B_T(s) = T(N_s * T).
\end{equation*}
\end{defn}

\begin{ex}
Suppose $M \subseteq \RR^d$ is a compact smooth submanifold, $dA$ is the
surface measure on $M$ and ${T_1} \phi = \int_M \phi \, dA$.  Then ${T_1}$
is a compactly supported distribution of order $0$.  If $\Re s > 0$,
then
\begin{equation*}
(N_s * {T_1})(u) = {T_1}(\tau_{-u}N_s) = \int_M N_s(v-u) \, dA(v),
\end{equation*}
where $\tau_{-u}$ is translation
and so
\begin{equation*}
B_{T_1}(s) = {T_1}(N_s * {T_1}) = \int_M \int_M N_s(v-u) \, dA(v) \, dA(u)
          = B_M(s),
\end{equation*}
where $B_M$ is the Brylinski beta function of $M$ as in
Definition 2.1 of \cite{bbs}.  Note that in this instance, the integrals
on the right hand sides make sense for $\Re s > -\dim(M)$.
\end{ex}

Now suppose $M \subseteq \RR^d$ is a compact smooth hypersurface, $dA$ is the
surface measure on $M$, and $\nu$ is the outward unit normal vector field
of $M$.
Define ${T_2} \phi = \int_M \partial_{\nu}\phi \, dA$.  Then $T_2$ is a
compactly supported distribution of order $1$.  If $\Re s > 1$, then
\begin{equation*}
(N_s * {T_2})(u) = {T_2}(\tau_{-u}N_s) = \int_M \partial_{\nu_v}N_s(v-u) \, dA(v),
\end{equation*}
and so
\begin{equation*}
B_{T_2}(s) = {T_2}(N_s * {T_2})
         = \int_M \int_M \partial_{\nu_u}\partial_{\nu_v}N_s(v-u) \, dA(v) \, dA(u)
\end{equation*}
Note that the integral
\begin{equation*}
\int_M \partial_{\nu_v}N_s(v-u) \, dA(v)
\end{equation*}
is defined for $\Re s > 2-d$.  In fact, it is defined even when $s=2-d$,
but is not continuous in $u$; in this case, it is proportional to the well
known {\em double layer potential} \cite{folland}.  Because of this
connection, we will refer to $B_{T_2}$ as the
{\em Brylinski beta function of a double layer on $M$}, and will use the
notation
\begin{equation*}
B_{\bar{\bar{M}}}=B_{T_2}.
\end{equation*}

\begin{rem}
In calculations, it is convenient to regard $N_s$ as a function of two
variables:
\begin{equation*}
N_s(u,v)=N_s(u-v).
\end{equation*}
\end{rem}

\begin{ex}
Let $M=S^{d-1}(R):=\{ u \in \RR^d\st \norm{u}=R \}$ be the sphere of
radius $R$ in $\RR^d$. Then 
\begin{equation*}
B_{\bar{\bar{S^{d-1}(R)}}}(s)
= \pi^{\frac{d-1}{2}}\omega_{d-1}\[
\frac{2^{s+d-5}R^{s+2d-4}s(s-2)\,(s+d-2)\Gamma\(\frac{s+d-3}{2}\)}
     {\Gamma\(\frac{s}{2}+d-1\)}
\]
\end{equation*}
where $\omega_n$ denotes the $n$-dimensional ``area" of $S^n(1)$.
\end{ex}

\begin{proof}

Suppose there is a uniform distribution of unit dipoles on $S^{d-1}(R)$ pointing
in the normal direction.  To compute its beta function, we need to integrate
the contribution from each pair of dipoles.  We fix one of the dipoles at
$v=(R,0,\dots,0)$, and integrate over the positions of the other dipole at
$u=(u_1,u_2,\dots,u_d) \in S^{d-1}(R)$.  By symmetry the total beta function is
simply $\omega_{d-1}R^{d-1}$ times the previous number.

Let $\nu_u$ and $\nu_v$ denote the normal vectors to $S^{d-1}(R)$ at
$u$ and $v$ respectively.  Then
\begin{equation*}
\begin{aligned}
\nu_u
=&\; \frac{u}{\norm{u}}=(u_1/R, u_2/R,\dots,u_d/R),
\qquad\text{and}\\
\nu_v
=&\; (1,0,0,\dots,0).
\end{aligned}
\end{equation*}

Therefore
\begin{equation*}
\begin{aligned}
(\partial_{\nu_u}\partial_{\nu_v}N_s)(u,(R,0,0,\dots,0))
=&\; (\nabla_u\partial_{\nu_v} N_s)(u,(R,0,0,\dots,0)) \cdot \nu_u\\
=&\; -s((u_1-R)^2+u_2^2+\cdots+u_d^2)^{\frac{s-2}{2}}\,{\frac{u_1}{R}}\\
 &\; +s(s-2)((u_1-R)^2+u_2^2+\cdots+u_d^2)^{\frac{s-4}{2}}(u_1-R)^2.
\end{aligned}
\end{equation*}
Using Euler's beta integral and the fact that
\begin{equation*}
\int_{-1}^1 x (1-x)^a (1-x^2)^{\frac{b}{2}} \, dx
= \frac{-2^{a+b+1} \,a\, \Gamma\(a+\frac b2+1\) \, \Gamma\(1+\frac b2\)}
          {\Gamma(a+b+3)},
\end{equation*}
we compute
\begin{equation*}
\begin{aligned}
&\;\int_{S^{d-1}(R)}(\partial_{\nu_u}\partial_{\nu_v}N_s)(u,(R,0,0,\dots,0)) \, dA(u)\\
= &\;\[
     \begin{gathered}
       \int_{-R}^R -s((u_1-R)^2+R^2-u_1^2)^{\frac{s-2}{2}}\,{\frac{u_1}{R}}
       \(\omega_{d-2}\,R\,\(R^2-u_1^2\)^{\frac{d-3}{2}}\, du_1\)\\                    
       +\int_{-R}^R s(s-2)((u_1-R)^2+R^2-u_1^2)^{\frac{s-4}{2}}(u_1-R)^2
                   \(\omega_{d-2}\,R\,\(R^2-u_1^2\)^{\frac{d-3}{2}}\, du_1\)
     \end{gathered}
     \]\\
= &\;2^{\frac{s-4}{2}} R^{s+d-3}\omega_{d-2} \[
\begin{gathered}
  s(s-2)\frac{2^{\frac{s}{2}+d-2}\Gamma(\frac{s+d-1}{2})\Gamma(\frac{d-1}{2})}
             {\Gamma(\frac{s}{2}+d-1)}
             + 2s \frac{\(\frac{s-2}{2}\) \, 2^{\frac{s}{2}+d-3} \,
                        \Gamma\(\frac{s+d-3}{2}\)\Gamma(\frac{d-1}{2})}
                       {\Gamma(\frac{s}{2}+d-1)}
\end{gathered}
\]\\
=&\; \pi^{\frac{d-1}{2}} \[
\frac{2^{s+d-4} R^{s+d-3} s(s-2)(s+d-2)\Gamma\(\frac{s+d-3}{2}\)}
     {\Gamma\(\frac{s}{2}+d-1\)}
\].\\
\end{aligned}
\end{equation*}
Therefore
\begin{equation*}
B_{\bar{\bar{S^{d-1}(R)}}}(s)
= \pi^{\frac{d-1}{2}}\omega_{d-1}\[
\frac{2^{s+d-5}R^{s+2d-4}s(s-2)\,(s+d-2)\Gamma\(\frac{s+d-3}{2}\)}
     {\Gamma\(\frac{s}{2}+d-1\)}
\].
\qedhere
\end{equation*}
\end{proof}

\section{The analytic continuation}\label{S:mainresult}

\begin{thm}\label{T:AC}
The function $B_{\bar{\bar{M}}}$ may be analytically continued to a meromorphic
function on $\CC$ with simple poles at $1-d-j$, $j=-2,0,2,4,\dots$.
Moreover, the residues are integrals over $M$
of polynomials in complete contractions of the covariant derivatives of
the second fundamental form.
\end{thm}

\begin{proof}
The proof is similar to the proof of \cite[Theorem 3.3]{bbs}, but the
details are more complicated; we give them here for completeness.
We find the analytic continuation of $B_{\bar{\bar{M}}}^v$; integration
over $M$ gives the stated result. 
If $\psi \in C_c^\infty(\RR^d)$ is identically one in a neighborhood
of $v$ then the localized beta function
\begin{equation*}
B_{\bar{\bar{M}}}^\psi(s) =
\int_{M} \partial_{\nu_u} \partial_{\nu_v} \norm{v-u}^s  \, \psi(u) \, dA(u)
\end{equation*}
has the same principal part as $B_{\bar{\bar{M}}}^v$ because their difference
extends to a holomorphic function on $\CC$.  So it suffices to prove
the result with $B_{\bar{\bar{M}}}^v$ replaced by $B_{\bar{\bar{M}}}^\psi$ for an
appropriate $\psi$.  By rotating and translating $M$, we may assume that $v=0$
and the tangent space to $M$ is $\RR^{d-1} \subseteq \RR^d$ (clearly this
process does not affect the beta functions).  Then, in a neighborhood
of $0$, $M$ is the graph of a function $f:\RR^{d-1} \to \RR^1$ which
vanishes to second order at $0$.  By making the neighborhood smaller,
we may assume $\norm{f(w)} < \norm{w}$.Choose $\psi$ to have support
in this neighborhood.  Then

\begin{equation*}
\begin{aligned}
\partial_{\nu_u}\partial_{\nu_v}N_s(u,0)
=&\; \frac{1}{\sqrt{1+\sum_{i=1}^{d-1}\(\frac{\partial f}{\partial u_i}\)^2}}
       \[
       \begin{gathered}
         -s(s-2)f(u_1,\dots,u_{d-1})
         \(\sum_{i=1}^{d-1}u_i^2+ f(u_1,\dots,u_{d-1})^2\)^{\frac{s-4}{2}}\\
         \(\sum_{i=1}^{d-1}{-\frac{\partial f}{\partial u_i}}u_i
           + f(u_1,\dots,u_{d-1})\)\\
           -s\(\sum_{i=1}^{d-1}u_i^2+ f(u_1,\dots,u_{d-1})^2\)^{\frac{s-2}{2}}
       \end{gathered} 
       \]
\end{aligned}
\end{equation*}
and
\begin{equation*}
  dA=\sqrt{1+\sum_{i=1}^{d-1}\(\frac{\partial f}{\partial u_i}\)^2}\,
    du_1\cdots du_{d-1}.
\end{equation*}
Therefore
\begin{equation*}
\begin{aligned}
 B_{\bar{\bar{M}}}^\psi(s)
=&\; \[
\begin{gathered}
  \int_{\RR^{d-1}} s(s-2) f(w)\(\norm{w}^2 + \norm{f(w)}^2\)^{\frac{s-4}{2}}
       \(-f(w)+\sum_{i=1}^{d-1}{\frac{\partial f}{\partial u_i}}u_i\)\,
       \phi(w) \, dw\\
       -\int_{\RR^{d-1}} s\(\norm{w}^2 + \norm{f(w)}^2\)^{\frac{s-2}{2}} \,
       \phi(w) \, dw
\end{gathered}
\]\\
\end{aligned}
\end{equation*}
where $\phi(w)=\psi(w, f(w))$.

Now,
\begin{equation*}
\begin{aligned}
B_{\bar{\bar{M}}}^\psi(s)
=&\; \[
     \begin{gathered}
       -s(s-2)\int_{\RR^{d-1}}\norm{w}^s
       \(1+{\frac{f(w)^2}{\norm{w}^2}}\)^{\frac{s-4}{2}}
       f(w)
       {\frac{(f(w)-\sum_{i=1}^{d-1}{\frac{\partial f}{\partial u_i}}u_i)}
             {\norm{w}^4}} \,\phi(w)\, dw\\
       -s \int_{\RR^{d-1}}\norm{w}^{s-2}
          \(1+{\frac{f(w)^2}{\norm{w}^2}}\)^{\frac{s-2}{2}} \,\phi(w)\, dw
     \end{gathered}
     \]\\
=&\; -s(s-2)\int_0^\infty r^{s+d-2} S_1(r,s)\,dr 
     - s\int_0^\infty r^{s+d-4} S_2(r,s)\,dr
\end{aligned}
\end{equation*}
where
\begin{equation*}
\begin{gathered}
S_1(r,s)=
\int_{S^{d-2}(1)}\(1+{\frac{f(rw)^2}{r^2}}\)^{\frac{s-4}{2}}
f(rw)
{\frac{(f(rw)-r\sum_{i=1}^{d-1}{\frac{\partial f}{\partial u_i}}(rw)\ u_i)}{r^4}}
\,\phi(w)\, d\sigma(w),\\
S_2(r,s)=
\int_{S^{d-2}(1)}\(1+{\frac{f(rw)^2}{r^2}}\)^{\frac{s-2}{2}}  \,\phi(w)\, d\sigma(w),
\end{gathered}
\end{equation*}
and $\sigma$ is the surface measure on $S^{d-2}(1)$.  Note that
$\frac{S_1(r,s)}{\omega_{d-2}}$ and $\frac{S_2(r,s)}{\omega_{d-2}}$ are the mean
values of the functions
$\(1 + \frac{\norm{f(rw)}^2}{r^2}\)^{\frac{s-4}{2}} \, \phi(w)$ and
$\(1 + \frac{\norm{f(rw)}^2}{r^2}\)^{\frac{s-2}{2}} \, \phi(w)$ on the sphere
of radius $r>0$. Define
$P: \RR \x S^{d-2}(1) \to \RR^{d-1}$ by $P(r, w) = rw$.
We can extend the definition of $S_1(r,s)$ and $S_2(r,s)$ to all real
values of $r$ by writing
\begin{equation*}
S_1(r,s) = \int_{S^{d-2}(1)} G_1(r, w, s) \, d\sigma(w)
\end{equation*}
and
\begin{equation*}
S_2(r,s) = \int_{S^{d-2}(1)} G_2(r, w, s) \, d\sigma(w),
\end{equation*}
where
\begin{equation*}
\begin{gathered}
G_1(r,w,s)
= \(1 + \frac{\norm{(f \circ P)(r,w)}^2}{r^2}\)^{\frac{s-4}{2}}
  {\frac{f \circ P(r,w)  g \circ P(r,w)}{r^4}}
  \, (\phi \circ P)(r, w),\\
g(rw) = f(rw)-r \sum_{i=1}^{n}{\frac{\partial f}{\partial u_i}}(rw)\ u_i,
\qquad\text{and}\\
G_2(r,w,s)
= \(1 + \frac{\norm{(f \circ P)(r,w)}^2}{r^2}\)^{\frac{s-2}{2}}
  \, (\phi \circ P)(r, w), \quad r \in \RR, \, s \in \CC.
\end{gathered}
\end{equation*}

By Lemma 3.2 of \cite{bbs}, $G_1,G_2 \in C^\infty(\RR \x S^{d-2}(1) \x \CC)$.
Moreover, for each $r \in \RR$ and $w \in S^{d-2}(1)$, the functions
$G_1(r, w,\cdot)$ and $G_2(r, w,\cdot)$ are entire, so by Lemma 3.1 of
\cite{bbs}, $S_1,S_2 \in C^\infty(\RR \x \CC)$, and $S_1(r, \cdot)$ and
$S_2(r, \cdot)$ are entire for each $r \in \RR$.  Note that by the equality
of mixed-partials, the functions
$\frac{\partial^j S_1}{\partial r^j}(r, \cdot)$ and
$\frac{\partial^j S_2}{\partial r^j}(r, \cdot)$ are also entire. Since
$G_1(-r, -w, s) = G_1(r,w,s)$ and $G_2(-r, -w, s) = G_2(r,w,s)$, it follows
that $S_1(\cdot, s)$ and $S_2(\cdot, s)$ are even, so
$\frac{\partial^j S_1}{\partial r^j}(0, \cdot) = 0$ and
$\frac{\partial^j S_2}{\partial r^j}(0, \cdot) = 0$ for odd $j$.

Now fix a positive integer $k$.  For $\Re s > 1-d$, we have
\begin{equation*}
\begin{aligned}
B_{\bar{\bar{M}}}^\psi(s)
= & \; -s(s-2) 
\int_0^1 r^{s+d-2} \[ S_1(r, s) - S_1(0,s)
       - r \frac{\partial S_1}{\partial r}(0,s) - \cdots
       - \frac{r^{k-1}}{(k-1)!} \frac{\partial^{k-1} S_1}{\partial r^{k-1}}(0,s) \]
         \, dr \\
   & \; -s(s-2) \int_1^\infty r^{s+d-2} S_1(r,s) \, dr
        -s(s-2) \sum_{j=0}^{k-1} \frac{1}{j! (s+d+j-1)}
          \frac{\partial^j S_1}{\partial r^j}(0,s) \\ 
  & \; -s 
\int_0^1 r^{s+d-4} \[ S_2(r, s) - S_2(0,s)
       - r \frac{\partial S_2}{\partial r}(0,s) - \cdots
       - \frac{r^{k+1}}{(k+1)!} \frac{\partial^{k+1} S_1}{\partial r^{k+1}}(0,s) \]
         \, dr \\
   & \; - s\int_1^\infty r^{s+d-4} S_2(r,s) \, dr
        -s \sum_{j=0}^{k+1} \frac{1}{j! (s+d+j-3)}
          \frac{\partial^j S_2}{\partial r^j}(0,s).               
\end{aligned}
\end{equation*}
By Taylor's theorem, the first integral on the right is defined and
holomorphic as a function of $s$ for $\Re s > 1-d-k$, so the right hand
side is a meromorphic function with only simple poles at $3-d$ and $1-d-j$
on the half-plane $\Re s > 1-d-k$.
Since
$\frac{\partial^j S_i}{\partial r^j}(0, \cdot) = 0$ ($i=1,2$), for odd $j$,
it follows that $B_{\bar{\bar{M}}}^\psi$ does not actually have
a pole at $1-d-j$ for odd $j$.  Since $k$ is arbitrary, this provides
the desired analytic continuation of $B_{\bar{\bar{M}}}^\psi$ to $\CC$.  Observe
that
\begin{equation*}
\begin{gathered}
 \Res_{s=3-d} B_{\bar{\bar{M}}}^v = -sS_2(0,s),
 \qquad\text{and}\\
\Res_{s=1-d-j} B_{\bar{\bar{M}}}^v
= -s(s-2)\frac{1}{j!} \frac{\partial^j S_1}{\partial r^j}(0, 1-d-j)
  -s \frac{1}{j+2!} \frac{\partial^{j+2} S_2}{\partial r^{j+2}}(0, 1-d-j).
\end{gathered}
\end{equation*}
  
Let $\mathrm{I}$ and $\mathrm{II}$ denote the first and second
fundamental forms of $M$.  Using the local parametrization $(w, f(w))$
of $M$, we find
\begin{equation*}
\mathrm{I}_{ij} = \delta_{ij} + \partial_i f \partial_j f, \quad\text{and}\quad
\mathrm{II}_{ij} = \frac{\partial_i \partial_j f}{(1+\norm{\nabla f}^2)^{1/2}}.
\end{equation*}
If $\alpha$ is a multi-index, let $\mathrm{II}_{ij;\alpha}$ denote the
$\alpha$-th covariant derivative of the second fundamental form.  Then
\begin{equation*}
(\mathrm{II}_{ij;\alpha})(0) = (\partial^\alpha \partial_i \partial_j f)(0) + Q
\end{equation*}
where $Q$ is a polynomial in the Taylor coefficients (at $0$) of $f$
of order less than or equal to $\abs{\alpha}+1$.  Such a ``triangular"
relation may be inverted to express the Taylor coefficients of $f$ as
polynomials in $(\mathrm{II}_{ij;\alpha})(0)$.

It follows that $\Res_{s=1-d-j} B_{\bar{\bar{M}}}^v$ may be expressed as a polynomial
in $(\mathrm{II}_{ij;\alpha})(0)$.
Let
\begin{equation*}
\M =
\left\{
M\subseteq \RR^d \; \middle|\;
\begin{aligned}
&\text{$M$ is an oriented hypersurface such that}\
  0\in M \ \text{and} \\
& T_0M = \{(x_1,x_2,\dots,x_{d-1},0)|\ x_i \in \RR\}
\end{aligned}
\right\}
\end{equation*}
Then $\mathrm{O}(d-2)$ acts on $\M$.  Moreover, the function
$\M\to\CC$ given by $M\mapsto \Res_{s=1-d-j} B_{\bar{\bar{M}}}^v$ is
invariant under this action.
By Weyl's theorem on invariants of the Orthogonal group, it
follows that this polynomial may be re-expressed as a polynomial in the
complete contractions of $\mathrm{II}_{ij;\alpha}(0)$
(cf. \cite{Gray,Weyl}).
\end{proof}

\section{Curves in $\RR^2$}\label{S:curves}

In this section, assume $M \subseteq \RR^2$ is a smooth curve.  We will
compute $\mathrm{Res}_{s=k} B_{\bar{\bar{M}}}$ for $k=1$, $-1$, $-3$
and $-5$ in terms of the curvature and its derivatives with respect to
a unit speed parametrization.  In this case, the second fundamental
form has only one component, namely curvature, and the covariant
derivatives of the second fundamental form are derivatives of
the curvature with respect to a unit speed parametrization.

\begin{thm}
\begin{equation*}
\begin{aligned}
\Res_{s=1} B_{\bar{\bar{M}}} &= -2 \, \mathrm{Length}(M)\\
\Res_{s=-1} B_{\bar{\bar{M}}} &= \int_M \frac{3}{4} \kappa_0^2  \, dA\\
\Res_{s=-3} B_{\bar{\bar{M}}} &=  \int_M \(
\frac{45}{64} \kappa_0^4 + \frac{5}{4} \kappa_1^2 
+\frac{15}{8} \kappa_0 \kappa_2 
\) \, dA\\
\Res_{s=-5} B_{\bar{\bar{M}}} &=  \int_M \(
\frac{175}{512} \kappa_0^6 + \frac{1295}{288} \kappa_0^2 \kappa_1^2 
+\frac{1085}{576} \kappa_0^3 \kappa_2 + \frac{175}{576} \kappa_2^2
+ \frac{35}{72} \kappa_1 \kappa_3
+ \frac{35}{144} \kappa_0 \kappa_4
\) \, dA
\end{aligned}
\end{equation*}
where $\kappa_0$ is the curvature and $\kappa_n$ is the $n$-th
derivative of $\kappa_0$ with respect to a unit speed parametrization.
\end{thm}

\begin{proof}
As in the proof of Theorem \ref{T:AC}, we compute the residues of
$B_{\bar{\bar{M}}}^v$ for each $v \in M$; integration over $M$ gives the stated
result. Also, we may assume that $v=0$, and locally $M$ is the graph of a
function $u_2=f(u_1)$ which vanishes to
second order at $0$.  By the proof of Theorem \ref{T:AC},
\begin{equation*}
\begin{gathered}
\Res_{s=1} B_{\bar{\bar{M}}}^v = -sS_2(0,s),
\qquad\text{and}\\
\Res_{s=-1-j} B_{\bar{\bar{M}}}^v
=-s(s-2)\frac{1}{j!} \frac{\partial^j S_1}{\partial r^j}(0, -1-j)
 -s\frac{1}{j+2!} \frac{\partial^{j+2} S_2}{\partial r^{j+2}}(0, -1-j),
\end{gathered}
\end{equation*}
where
\begin{equation*}
\begin{gathered}
S_1(r,s)
=\sum_{w=\pm 1}\(1+{\frac{f(rw)^2}{r^2}}\)^{\frac{s-4}{2}}
 f(rw){\frac{f(rw)-rwf'(rw)}{r^4}} \,\phi(rw),
\qquad\text{and}\\
S_2(r,s)=\sum_{w=\pm 1}\(1+{\frac{f(rw)^2}{r^2}}\)^{\frac{s-2}{2}}  \,\phi(rw)
\end{gathered}
\end{equation*}
for small $r$.

If we write
\begin{equation*}
\begin{gathered}
f(u_1) = \\
a_2 u_1^2 + a_3 u_1^3 + a_4 u_1^4 + \\
a_5 u_1^5 + a_6 u_1^6 + a_7 u_1^7 + a_8 u_1^8 
+ O(u_1^9),\\
\(1+\frac{f(rw)^2}{r^2}\)^{\frac{s-4}{2}}
f(rw)
\frac{f(rw)-rw f'(rw)}{r^4}\\
=
\(1 + \frac{(s-4)}{2} \frac{f^2}{r^2} + \frac{(s-4)(s-6)}{8} \frac{f^4}{r^4}+
\frac{(s-4)(s-6)(s-8)}{48} \frac{f^6}{r^6}\)
f(rw)
\frac{f(rw)-rwf'(rw)}{r^4}+ O(r^8),
\end{gathered}
\end{equation*}
and
\begin{equation*}
\begin{gathered}
\(1+{\frac{f(rw)^2}{r^2}}\)^{\frac{s-2}{2}}\\
=1 +  \frac{(s-2)}{2} \frac{f^2}{r^2} + \frac{(s-2)(s-4)}{8} \frac{f^4}{r^4} +
\frac{(s-2)(s-4)(s-6)}{48} \frac{f^6}{r^6}+ O(r^8),
\end{gathered}
\end{equation*}
substitute into the definition of $S_1(r,s)$ and $S_2(r,s)$ and perform the
indicated differentiation (using the computer algebra package
Maxima \cite{maxima}), we find
\begin{equation}\label{E:residues-basis-curves}
\begin{aligned}
\Res_{s=1} B_{\bar{\bar{M}}}^v
=&\; -2\\
\Res_{s=-1} B_{\bar{\bar{M}}}^v
=&\; 3\,a_{2}^2\\
\Res_{s=-3} B_{\bar{\bar{M}}}^v
=&\; {{360\,a_{2}\,a_{4}+180\,a_{3}^2-315\,a_{2}^4}\over{4}}\\
\Res_{s=-5} B_{\bar{\bar{M}}}^v
=&\;{{2800\,a_{2}\,a_{6}+2800\,a_{3}\,a_{5}+1400\,a_{4}^2
    -12600\,a_{2}^3\,a_{4}-18900\,a_{2}^2\,a_{3}^2+5775\,a_{2}^6}\over{8}}
\end{aligned}
\end{equation}

At $u_1=0$, we compute (using Maxima)
\begin{equation*}
\begin{aligned}
\kappa_0
=&\; 2\,a_{2}\\
\kappa_1
=&\; 6\,a_{3}\\
\kappa_2
=&\; 24\,a_{4}-24\,a_{2}^3\\
\kappa_3
=&\; 120\,a_{5}-456\,a_{2}^2\,a_{3}\\
\kappa_4
=&\; 720\,a_{6}-16\,a_{2}^2\,\left(24\,a_{4}-24\,a_{2}^3\right)
     -2880\, a_{2}^2\,a_{4}-3456\,a_{2}\,a_{3}^2+1440\,a_{2}^5\\
\kappa_5
=&\; 5040\,a_{7}-24\,a_{2}^2\,\left(120\,a_{5}-432\,a_{2}^2\,a_{3}\right)
     -16\,a_{2}^2\,\left(120\,a_{5}-456\,a_{2}^2\,a_{3}\right)
     -21600\,a_{2}^2\,a_{5}\\
 &\; -180\,a_{2}\,a_{3}\,\left(24\,a_{4}-24\,a_{2}^3 \right)
     -52992\,a_{2}\,a_{3}\,a_{4}-10368\,a_{3}^3+65664\,a_{2}^4\, a_{3}.
     \qedhere
\end{aligned}
\end{equation*}
Solving for $a_k$ in terms of $\kappa_0, \dots ,\kappa_{k-2}$ and substituting
in equation (\ref{E:residues-basis-curves}) gives the stated result.
\end{proof}

\section{Surfaces in $\RR^3$}\label{S:surfaces}

In this section, assume $M \subseteq \RR^3$ is a smooth surface.  We will
compute $\mathrm{Res}_{s=k} B_{\bar{\bar{M}}}$ for $k=0$, $-2$ and $-4$ in terms
of $\I$, $\II$ and the first two covariant derivatives of
$\mathrm{II}$.

Fix $v \in M$ and choose coordinates such that $\I_{ij} = \delta_{ij}$
at $v$.  Let $\II_{ij;k}$ and $\II_{ij;kl}$ denote the components of
the first and second covariant derivatives of $\II$.  Define, at $v$
\begin{equation*}
\begin{aligned}
H_0 &= \II_{ii} \\
H_1 &= \II_{ij} \II_{ij} \\
H_2 &= \II_{ij;k} \II_{ij;k} \\
H_3 &= \II_{ii;jj} \\
H_4 &= \II_{ij} \II_{ij;kk} \\
H_5 &= \II_{ij} \II_{kk;ij} \\
H_6 &= \II_{ii;k} \II_{jj;k}.
\end{aligned}
\end{equation*}
Here we are using the extended Einstein summation convention, where we
sum over repeated indices, even if they are both covariant.  This is
justified because we are working in a coordinate system which is
orthogonal at $v$.  One can get formulae for the $H_\alpha$ in a
general coordinate system by first raising one of the repeated indices.

Being complete contractions of tensors, $H_\alpha$, $\alpha=0, \dots, 6$
are smooth functions on $M$, independent of any coordinate choices.
Note that $H_0$ is just the mean-curvature.

\begin{thm}
\begin{equation*}
\begin{aligned}
\Res_{s=0} B_{\bar{\bar{M}}}
=&\; 0\\
\Res_{s=-2} B_{\bar{\bar{M}}}
=&\; \frac{\pi}{4} \int_M (H_0^2 + 2H_1) \, dA\\
\Res_{s=-4} B_{\bar{\bar{M}}}
=&\; \frac{\pi}{8} \int_M \(
-\frac{15}{16} H_0^4 + \frac{15}{4} H_1^2 + 2 H_2
+\frac{3}{2} H_0 H_3 + \frac{9}{2} H_4 + \frac{3}{2} H_5
+3 H_6
\) \, dA
\end{aligned}
\end{equation*}
\end{thm}

\begin{proof}
As in the proof of Theorem \ref{T:AC}, we compute the residues of
$B_{\bar{\bar{M}}}^v$ for each $v \in M$; integration over $M$ gives the stated
result.  Assume that $v=0$ and the $u_1u_2$-plane is tangent to $M$ at $0$.
So locally $M$ is the graph of a function $u_3=f(w)$, where $w=(u_1,u_2)$,
which vanishes to second order at $0$.  By the proof of Theorem \ref{T:AC},
\begin{equation*}
\begin{gathered}
\Res_{s=0} B_{\bar{\bar{M}}}^v = -sS_2(0,s),
\qquad\text{and}\\
\Res_{s=-2-j} B_{\bar{\bar{M}}}^v
= -s(s-2)\frac{1}{j!} \frac{\partial^j S_1}{\partial r^j}(0, -2-j)
  -s\frac{1}{(j+2)!} \frac{\partial^{j+2} S_2}{\partial r^{j+2}}(0, -2-j),
\end{gathered}
\end{equation*}
where
\begin{equation*}
S_1(r,s)
=\int_{S^{1}(1)}\(1+{\frac{f(rw)^2}{r^2}}\)^{\frac{s-4}{2}}
 f(rw)
 {\frac{(f(rw)-r\sum_{i=1}^{2}{\frac{\partial f}{\partial u_i}}(rw)\ u_i)}{r^4}}
 \,\phi(w)\,d\sigma(w),\\
\end{equation*}
and
\begin{equation*}
S_2(r,s)
=\int_{S^{1}(1)}\(1+{\frac{f(rw)^2}{r^2}}\)^{\frac{s-2}{2}}  \,\phi(w)\, d\sigma(w)
\end{equation*}
for small $r$.

If we write
\begin{equation}\label{E:monge}
\begin{gathered}
f(u_1, u_2) = \\
b_1 u_1^2 + b_2 u_1 u_2 + b_3 u_2^2 + \\
c_1 u_1^3 + c_2 u_1^2 u_2 + c_3 u_1 u_2^2 + c_4 u_2^3 + \\
d_1 u_1^4 + d_2 u_1^3 u_2 + d_3 u_1^2 u_2^2 + d_4 u_1 u_2^3 + d_5 u_2^4
+ O(\norm{w}^5),\\
\end{gathered}
\end{equation}
\begin{equation*}
\begin{gathered}
\(1+{\frac{f(rw)^2}{r^2}}\)^{\frac{s-4}{2}}
f(rw)
{\frac{(f(rw)-r \sum_{i=1}^{2}{\frac{\partial f}{\partial u_i}}(rw)\ u_i)}{r^4}}\\
=\(1 +  \frac{(s-4)}{2} \frac{f^2}{r^2} + \frac{(s-4)(s-6)}{8} \frac{f^4}{r^4}\)
 f
 {\frac{(f-r \sum_{i=1}^{2}{\frac{\partial f}{\partial u_i}}(rw) \ u_i)}{r^4}}
 + O(r^6),
\end{gathered}
\end{equation*}
and
\begin{equation*}
\begin{gathered}
\(1+{\frac{f(rw)^2}{r^2}}\)^{\frac{s-2}{2}}\\
=1 +  \frac{(s-2)}{2} \frac{f^2}{r^2} + \frac{(s-2)(s-4)}{8} \frac{f^4}{r^4}+
O(r^6),
\end{gathered}
\end{equation*}
substitute into the definition of $S_1(r,s)$ and $S_2(r,s)$ and perform the
indicated differentiation (using Maxima), we find
\begin{equation}\label{E:residues-basis}
\begin{aligned}
\Res_{s=0} B_{\bar{\bar{M}}}^v
=&\; 0\\
\Res_{s=-2} B_{\bar{\bar{M}}}^v
=&\; \pi (3b_3^2 + 2 b_1 b_3 + b_2^2 + 3b_1^2)\\
\Res_{s=-4} B_{\bar{\bar{M}}}^v 
=&\; \frac{\pi}{8}
\(
\begin{gathered}
360\,b_{3}\,d_{5}+72\,b_{1}\,d_{5}+72\,b_{2}\,d_{4}+180\,c_{4}^2\\
+72\,c_{2}\,c_{4}+72\,b_{3}\,d_{3} +72\,b_{1}\,d_{3}+36\,c_{3}^2\\
+72\,c_{1}\,c_{3}-315\,b_{3}^4-180\,b_{1}\,b_{3}^3-270\,b_{2}^2\,b_{3}^2\\
-162 \,b_{1}^2\,b_{3}^2-324\,b_{1}\,b_{2}^2\,b_{3}+72\,d_{1}\,b_{3}\\
-180\, b_{1}^3\,b_{3}+72\,b_{2}\,d_{2}+36\,c_{2}^2-27\,b_{2}^4\\
-270\,b_{1}^2\,b_{2}^2+360\,b_{1}\,d_{1}+180\,c_{1}^2-315\,b_{1}^4
\end{gathered}
\)
\end{aligned}
\end{equation}

In order to re-express these in terms of complete contractions of the
covariant derivatives of the second fundamental form, we introduce the
concept of {\em weight}.

Observe that if $\lambda\in\RR^+$ and $M\in \M$, the set
\begin{equation*}
\lambda^{-1}M=\{\lambda^{-1}u \st u \in M\}
\end{equation*}
also belongs to $\M$, and moreover the map
\begin{equation*}
(\lambda, M) \mapsto \lambda^{-1}M
\end{equation*}
defines an action of the group $\RR^+$ on $\M$.

Let $\F$ be a polynomial in the components of the covariant derivatives at
$0$ of the second fundamental form of $M$ (in the $u_1u_2$-coordinate system).
We will regard $\F$ as a map $\F:\M\to\CC$.

If $\F$ satisfies $\F(\lambda^{-1}M)=\lambda^l\F(M)$
for all $\lambda\in\RR^+$ and $M\in\M$, we will say that $\F$ is
a functional of
{\em weight} $l$.  Observe that functionals of weight $l$ form a
finite dimensional complex vector space, and moreover, if $\F$ has
weight $l$ and $\G$ has weight $m$ then $\F\G$ has weight $l+m$.

Now let us consider the weights of the various functionals at hand.

\begin{enumerate}

\item
The components $\mathrm{II}_{ij;\alpha}(0)$ (in the $u_1u_2$-coordinate system)
have weight $\abs{\alpha}+1$.

\item
When we express $M$ locally as the graph of a function as in
equation (\ref{E:monge}), the monomial $b_i^\alpha c_j^\beta d_k^\gamma$
has weight $\alpha+2\beta+3\gamma$, and the set of all such monomials
forms a basis for the space of functionals of weight
$\alpha+2\beta+3\gamma$.


\item
The functional $M\mapsto\Res_{s=-k} B_{\bar{\bar{M}}}^v$
has weight $k$.

\end{enumerate}

Observe that every monomial in complete contractions of
$\mathrm{II}_{ij;\alpha}(0)$ has a weight.  According to Weyl's theorem,
the vector space of $\mathrm{O}(2)$-invariant functionals of weight $l$
is spanned by the weight $l$ monomials in complete contractions of
$\mathrm{II}_{ij;\alpha}(0)$.  For small values of the weight, it is
easy to enumerate them:



\begin{center}
  \begin{tabular}{||c|c||}
    \hline
    \hline
    Weight & Monomial\\
    \hline
    $0$    & $1$\\
    \hline
    $1$    & $H_0$\\
    \hline
    $2$    & $H_0^2, H_1$\\
    \hline
    $3$    & $H_0^3, H_0 H_1, H_3$\\
    \hline
    $4$    & $H_0^4, H_1^2, H_0^2H_1, H_2, H_0H_3, H_4, H_5, H_6$\\
    \hline
    \hline
  \end{tabular}
\end{center}

From this we conclude that $\Res_{s=0} B_{\bar{\bar{M}}}^v$ is independent
of $M$, $\Res_{s=-2} B_{\bar{\bar{M}}}^v$ is a linear combination of
$H_0^2$ and $H_1$, and $\Res_{s=-4} B_{\bar{\bar{M}}}^v$ is a linear combination of
$H_0^4$, $ H_1^2$, $ H_0^2H_1$, $ H_2$, $ H_0H_3$, $ H_4$, $ H_5$, $ H_6$

Using Maxima, we calculate the components of $\I$ and $\II$ with respect
to the coordinates $u_1, u_2$ on $M$ to order two near $0$, and so determine
the first and second covariant derivatives of $\II$ at $0$.
Using this,
we find that, at $0$, we have
\begin{equation*}
\begin{aligned}
H_0^2
=&\;4b_3^2+8b_1b_3+4b_1^2 \\
H_1
=&\;4b_3^2+2b_2^2+4b_1^2 \\
H_0^4
=&\;16b_3^4+64b_1b_3^3+96b_1^2b_3^2+64b_1^3b_3+16b_1^4 \\
H_1^2
=&\;16b_3^4+16b_2^2b_3^2+32b_1^2b_3^2+4b_2^4+
    16b_1^2b_2^2+16b_1^4 \\
H_0^2H_1
=&\;16b_3^4+32b_1b_3^3+8b_2^2b_3^2+32b_1^2b_3^2+
     16b_1b_2^2b_3+32b_1^3b_3+8b_1^2b_2^2+16b_1^4 \\
H_2
=&\;36c_4^2+12c_3^2+12c_2^2+36c_1^2 \\
H_0H_3
=&\; 48b_3d_5+48b_1d_5+16b_3d_3+16b_1d_3-48b_3^4\\
 &\; -64b_1b_3^3-32b_2^2b_3^2-32b_1^2b_3^2-64b_1b_2^2b_3\\
 &\; +48d_1b_3-64b_1^3b_3-32b_1^2b_2^2+48b_1d_1-48b_1^4\\
H_4
=&\; 48b_3d_5+12b_2d_4+8b_3d_3+8b_1d_3-48b_3^4\\
 &\; -48b_2^2b_3^2-32b_1^2b_3^2-32b_1b_2^2b_3+12b_2d_2\\
 &\; -8b_2^4-48b_1^2b_2^2+48b_1d_1-48b_1^4 \\
H_5
=&\; 48b_3d_5+12b_2d_4+8b_3d_3+8b_1d_3-48b_3^4\\
 &\; -16b_1b_3^3-44b_2^2b_3^2-56b_1b_2^2b_3-16b_1^3b_3\\
 &\; +12b_2d_2-4b_2^4-44b_1^2b_2^2+48b_1d_1-48b_1^4 \\
H_6
=&\;36c_4^2+24c_2c_4+4c_3^2+24c_1c_3+4c_2^2+36c_1^2
\end{aligned}
\end{equation*}

Writing
\begin{equation*}
\begin{aligned}
\Res_{s=-2} B_{\bar{\bar{M}}}^v
=&\; \beta_1 H_0^2 +  \beta_2 H_1,\qquad\text{and}\\
\Res_{s=-4} B_{\bar{\bar{M}}}^v
=&\;
\gamma_1 H_0^4 + \gamma_2 H_1^2 + \gamma_3 H_0^2H_1 + \gamma_4 H_2
+ \gamma_5 H_0H_3 + \gamma_6 H_4 + \gamma_7 H_5 + \gamma_8 H_6
\end{aligned}
\end{equation*}
and comparing with equation (\ref{E:residues-basis}) gives an over-determined
system of linear equations for the $\beta_i$, $i=1,2$, and $\gamma_i$,
$i=1,\dots, 8$, which has a solution by Weyl's theorem.  Solving them gives
the coefficients in the statement of the theorem.
\end{proof}

\bibliographystyle{amsplain}
\bibliography{v9-bbdl}

\end{document}